\documentclass[a4paper,12pt]{article}
\usepackage[utf8]{inputenc}
\usepackage[]{amsmath}
\usepackage[]{amssymb}
\usepackage[]{amsthm}
\usepackage{stmaryrd}
\usepackage{graphicx}
\author{Magdalena Zielenkiewicz}
\title{Residue formulas for push-forwards in equivariant cohomology - a symplectic approach}
\date{}
\newcommand{\git}{/\mkern-6mu/}
\newcommand{\Z}{\mathbb{Z}}
\newcommand{\R}{\mathbb{R}}
\newcommand{\C}{\mathbb{C}}

\def\T{S}
\def\C{\mathbb{C}}
\def\R{\mathcal{R}}
\def\RR{\mathbb{R}}

\def\P{\mathbb{P}}
\def\E{\mathbb{E}}
\def\B{\mathbb{B}}
\def\X{X}

\def\lie#1{\mathfrak{#1}}

\def\MG{\X \git G}
\DeclareMathOperator{\vol}{vol}
\DeclareMathOperator{\Res}{Res}
\newtheorem{theorem}{Theorem}[section]

\newtheorem{proposition}[theorem]{Proposition}
\newtheorem{corollary}[theorem]{Corollary}

\theoremstyle{definition}
\newtheorem{remark}[theorem]{Remark}
\newtheorem{example}[theorem]{Example}
\numberwithin{equation}{section}

\begin{document}
\maketitle

\begin{abstract}
In \cite{guillemin1996} Guillemin and Kalkman proved how the nonabelian localization theorem of Jeffrey and Kirwan (\cite{jeffrey1995}) can be rephrased in terms of certain iterated residue maps, in the case of torus actions. In \cite{zielenkiewicz2014} we describe the push-forward in equivariant cohomology of homogeneous spaces of classical Lie groups, with the action of the maximal torus, in terms of iterated residues at infinity of certain complex variable functions. The aim of this paper is to show how, in the special case of classical Grassmannians, the residue formulas obtained in \cite{zielenkiewicz2014} can be deduced from the ones described in \cite{jeffrey1995} and \cite{guillemin1996}. 
\end{abstract}

\section{Introduction}

If $\X$ is a compact manifold with an action of a torus $\T$ and the fixed point set $\X^{\T}$ of the action is finite, the Atiyah--Bott--Berline--Vergne formula expresses the push-forward in equivariant cohomology of an element $\alpha \in H^*_\T(\X)$ as a finite sum of local contributions:
 \[\int_{\X} \alpha  = \sum_{p \in \X^\T} \frac{i_p^*\alpha}{e_p},\]
where $i_p: \{ p\} \to \X$ is the inclusion of the fixed point  and $e_p$ is the equivariant Euler class of the tangent bundle at $p$. \newline

In the case of the complex Grassmannian $Grass_m(\C^n)$ the fixed points are indexed by the partitions $\lambda = (\lambda_1 \geq \dots \geq \lambda_m)$ with $\lambda_i \in \Z_{> 0}$. Let us denote by $\alpha(\R)$ a characteristic class of the tautological bundle $\R$, then the above formula takes the form:
\[ \int_{Grass_{m}{\C^n}} \alpha(\R) = \sum_{p_\lambda} \frac{\alpha_{|_{p_\lambda}}}{e_{p_\lambda}} .\]
In \cite{zielenkiewicz2014} we have proven that if $\alpha(\R)$ restricted to the fixed points of the action is given by a symmetric polynomial $V$, this expression is equal to the following residue

\[ \int_{Grass_{m}{\C^n}} \alpha(\R) = \frac{1}{m!} Res_{\mathbf{z} = \infty} \frac{V(z_1,...,z_m)\prod_{i \neq j}(z_i - z_j)}{\prod_{i=1}^n\prod_{j=1}^m(t_i - z_j)}. \]
Similar results have been obtained for other types of Grassmannians (Lagrangian and orthogonal ones). The formulas are given in chapter \ref{wzorki}. In their paper \cite{berczi2012}, B\'erczi and Szenes found a formula for an integral over flag variety (\cite{berczi2012}, Chapter 6.3). All the obtained formulas (for classical, Lagrangian and orthogonal Grassmannians and the one derived by B\'erczi and Szenes) can be uniformly written as special cases of one formula, involving the action of the Weyl group on the characters of the natural representation of the torus action. The result is presented in \cite{zielenkiewicz2014}. \newline

These results, however, have been lacking a geometric motivation. In this paper we relate them to the nonabelian localization theorem by Jeffrey and Kirwan. In particular, we show how the residue-type formula for the classical Grassmannian can be obtained from generalization of the Jeffrey--Kirwan theorem, applied to the group $U(k)$ acting on the space $Hom(\C^k, \C^n)$. \newline

Chapter~\ref{notation} describes the notation conventions used in this paper, trying to find a compromise between the notations of the two influential papers \cite{jeffrey1995} and \cite{guillemin1996}. Chapter \ref{prelim} briefly describes localization theorems in equivariant cohomology, the results obtained in \cite{zielenkiewicz2014} and the various approaches to nonabelian localization. In Chapter~\ref{gk_equiv} we adapt the Jeffrey--Kirwan theorem to the $\T$-equivariant setting, for a torus $\T$, by modifying the proof of the Jeffrey--Kirwan theorem given by Guillemin and Kalkman in \cite{guillemin1996}. Two examples at the end of Chapter~\ref{gk_equiv} relate the equivariant version of Jeffrey--Kirwan theorem to the formulas in \cite{zielenkiewicz2014}. The Appendix shows the details of the adaptation of the Guillemin--Kalkman proof of the nonabelian localization theorem to the case of $S^1$-action. 

\subsection*{Acknowledgements}
This work was partially supported by NCN grant 2015/17/N/ST1/02327.

\section{Notation}\label{notation}
Almost every paper we refer to uses different notation conventions. We will use the following notation:

\begin{itemize}
	\item $\X$ is a symplectic manifold, usually compact 
	\item $G$ is a compact Lie group,
	\item $T$ is a maximal torus in $G$,
	\item $\mu_G: \X \to \mathfrak{g}$ is the moment map for the action of $G$ on $\X$,
	\item $\X \git G = \mu_{G}^{-1}(0)/G$ is the symplectic reduction, compact
	\item $\T$ is a torus (not related to $G$),
	\item $\kappa$ is a Kirwan map,
	\item $\kappa^{\T}$ is the $\T$-equivariant Kirwan map.
\end{itemize}

\section{Preliminaries}\label{prelim}

\subsection{Localization in equivariant cohomology}

Let $\X$ be a compact space equipped with an action of the torus $\T=(S^1)^n$. Let us consider the $\T$-equivariant cohomology of $\X$, $H^*_{\T}(\X)$, defined as 
\[ H^*_{\T}(\X) = H^*(\E\T \times^{\T} \X),\]
where $\E\T$ is a contractible space on which $\T$ acts freely (the total space of the universal $\T$-bundle). The space $\E\T$ is infinite dimensional, so one constructs finite dimensional approximations $\E_m$ such that for $m$ sufficiently large with respect to $i$ we have isomorphisms
\[H^i_{\T}(\X) \simeq H^i(\E_m \times^{\T} \X) \textrm{ for all } i \ll m. \]
We consider cohomology theories with coefficients in a field (usually $\C$). \newline

The characters of the torus action can be identified with elements of $H^*_{\T}(\X)$ as follows. For a character $\chi \in \T^{\#}$ let $\C_{\chi}$ denote the one-dimensional representation of $\T$ determined by $\chi$. Then $L_{\chi}: \E\T \times^{\T} \C_{\chi} \to \B\T $
is an equivariant line bundle, and to the character $\chi$ we can associate the first equivariant Chern class of the associated line bundle 
\[\chi \mapsto c_1^{\T}(L_{\chi}) \in H^2_{\T}(pt).\] 
We will abuse notation and use this assignment to identify characters with elements in equivariant cohomology.  
\newline

It turns out that if the set of fixed points of the action is finite, then after localizing with respect to the multiplicative set generated by the nonzero characters of the action the equivariant cohomology of $\X$ is isomorphic to the equivariant cohomology of the fixed point set $\X^{\T}$.

\begin{theorem}[Quillen (\cite{quillen1971})] Let $\X$ be a compact $\T$-space. The inclusion  $i: \X^\T \hookrightarrow \X$ induces an isomorphism
\[ i^*: H_{\T}^*(\X)[(\T^{\#} \setminus\{ 0\})^{-1}] \stackrel{\simeq}{\longrightarrow} H_{\T}^*(\X^\T)[(\T^{\#} \setminus\{ 0\})^{-1}]\]
after localizing with respect to the multiplicative system consisting of finite products of elements $c_1^{\T}(L_{\chi})$, for $\chi \in \T^{\#}\setminus \{0\}$.
\end{theorem}

For a proper map $f: \X \to Y$, one can define the equivariant Gysin homomorphism 
\[f_*^{\T} : H^i_{\T}(\X) \to H^{i+2d}_{\T}(Y), \]
where $d = \dim Y - \dim \X$. One defines $f_*^{\T}$ to be the standard non-equivariant Gysin homomorphism for the map $\E_m \times^{\T} \X \to \E_m \times^{\T} Y$. Gysin homomorphisms are also called push-forwards in cohomology, and if $f$ is a projection in a fiber bundle they can be interpreted as integration along fibers. The equivariant Gysin homomorphism associated with the map $\X \to pt$ is often denoted by $a \mapsto \int_{\X} a$. \newline

If $\X$ is a compact manifold and the fixed point set is finite, the localization theorem can be rephrased as follows.

\begin{theorem}[Atiyah--Bott (\cite{atiyah1984}), Berline--Vergne (\cite{berline1982})] Suppose $\X$ is a compact manifold with an action of a torus $\T$ such that $\# \X^\T < \infty$. Then for $\alpha \in H_{\T}^*(\X)$ one has
\[ \int_{\X} \alpha  = \sum_{p \in \X^\T} \frac{i_p^*\alpha}{e_p},\]
where $e_p$ is the equivariant Euler class of the tangent bundle at the fixed point $p$, and $i_p^*$ is induced by the inclusion of the fixed point $p$ into $\X$. 
\end{theorem}

To effectively compute the right-hand side of the above formula, one needs to know the equivariant Euler class at the fixed points of the action. The Euler class at the fixed point $p$ is given by the product of the weights of the torus action on the tangent space. If $\X$ is a homogeneous space $G/P$, where $P$ is a parabolic subgroup, the tangent space at $1$ is canonically isomorphic to $\lie{g}/\lie{p}$. The weights of the torus action are then the positive roots $\Phi^+ \setminus \Phi_P^+$.

\subsection{Push-forward in equivariant cohomology and residue formulas}\label{wzorki}

The Atiyah--Bott--Berline--Vergne formula expresses the push-forward in $\T$-equivariant cohomology as a finite sum of rational functions, indexed by the fixed points of the action. If $\X=G/P$ is a homogeneous space of a compact Lie group for a parabolic subgroup $P$, the fixed points are indexed by the quotient $W/W_P$ of the Weyl group $W$ of $G$ by the Weyl group $W_P$ of the parabolic subgroup $P$. In \cite{berczi2012}, B\'erczi and Szenes used iterated residues at infinity to compute the push-forward of a certain characteristic class over the complete flag variety. In \cite{zielenkiewicz2014} we prove similar residue-type formulas for push-forwards over Grassmannians (classical, Lagrangian and orthogonal). For example, the push-forward of a characteristic class $\alpha$ of the tautological bundle $\mathcal{R}$ over the Grassmannian $Grass_{m}(\C^n)$, which at the fixed points of the action is given by the symmetric polynomial $V$, can be expressed as

\[\int_{Grass_{m}(\C^n)} \!\!\!\! \alpha(\mathcal{R}) = \Res_{z_1 = \infty}\Res_{z_2 = \infty}...\Res_{z_m = \infty} \frac{1}{m!}\frac{V(z_1,...,z_m)\prod_{i \neq j}(z_i - z_j)}{\prod_{i=1}^n\prod_{j=1}^m(t_i - z_j)},\]
where $t_1,\dots,t_n$ are the characters of the action of the maximal torus, and $z_1, \dots, z_m$ are formal variables (whose interpretation and geometric meaning can be given using the Jeffrey--Kirwan nonabelian localization theorem, described in the next section). \newline

Similar formulas can be obtained for the Lagrangian Grassmannians $LG(n)$ and the orthogonal Grassmannians $OG(n,2n)$ and $OG(n, 2n+1)$, as follows. Denote by $\mathbf{z} = \{ z_1, z_2, \dots, z_n \}$. Then 

\[\int_{LG(n)} \alpha(\mathcal{R}) = \frac{1}{n!} \Res_{\mathbf{z} = \infty} \frac{V(z_1,...,z_n)\prod_{i<j}(z_j - z_i)}{\prod_{i=1}^n(t_i - z_i)(t_i + z_i)\prod_{i<j}(t_i + t_j)(t_j - t_i)}\]

\[ \int_{OG(n,2n)} \!\!\! \alpha(\mathcal{R}) = \frac{1}{n!} \Res_{\mathbf{z} = \infty} \frac{V(z_1,...,z_n)\prod_{i<j}(z_j - z_i)2^n z_1...z_n}{\prod_{i=1}^n(t_i - z_i)(t_i + z_i)\prod_{i<j}(t_i + t_j)(t_j - t_i)}\]

\[ \int_{OG(n,2n+1)} \!\!\! \alpha(\mathcal{R}) = \frac{1}{n!} \Res_{\mathbf{z} = \infty} \frac{V(z_1,...,z_n)\prod_{i<j}(z_j - z_i)2^n}{\prod_{i=1}^n(t_i - z_i)(t_i + z_i)\prod_{i<j}(t_i + t_j)(t_j - t_i)}. \]
These formulas can be rewritten in the following way, which is far less efficient for computations (due to the higher degrees of the polynomials appearing in the denominators), but has other advantages, described below.

\[\int_{LG(n)} \alpha(\mathcal{R}) = \Res_{\mathbf{z} = \infty} \frac{1}{n!}\frac{V(z_1,...,z_n)\prod_{i \neq j}(z_j - z_i)\prod_{i < j}(z_i+z_j)}{\prod_{i=1}^n\prod_{j=1}^n(t_i - z_j)(t_i + z_j)}.\]

\[ \int_{OG(n,2n)} \!\!\! \alpha(\mathcal{R}) = \Res_{\mathbf{z} = \infty} \frac{1}{n!}\frac{V(z_1,...,z_n)\prod_{i \neq j}(z_j - z_i)\prod_{i < j}(z_i+z_j)2^n\prod_{i=1}^n z_i}{\prod_{i=1}^n\prod_{j=1}^n(t_i - z_j)(t_i + z_j)}\]

\[ \int_{OG(n,2n+1)} \!\!\! \alpha(\mathcal{R}) = \Res_{\mathbf{z} = \infty} \frac{1}{n!}\frac{V(z_1,...,z_n)\prod_{i \neq j}(z_j - z_i)\prod_{i < j}(z_i+z_j)2^n}{\prod_{i=1}^n\prod_{j=1}^n(t_i - z_j)(t_i + z_j)}.\]

The advantage of these formulas, at first glance more complicated, is that one can now uniformly describe all of them, using one expression (which also agrees with the case of flag varieties investigated by B\'erczi and Szenes and may be true for arbitrary homogeneous spaces of Lie groups),

\[\int_{G/P} \alpha(\mathcal{R}) =  \Res_{\mathbf{z} = \infty} \frac{1}{|W_P|}\frac{V(z_1,...,z_n)\prod_{i=1}^n \prod_{x \in X_i \setminus \{z_i\} } (z_i - x)}{\prod_{i=1}^n\prod_{x \in X_i} (t_i - x ) \prod_{y \in \Phi^+ \setminus {\Phi_P}^+} y }.\]
Here we push-forward a characteristic class of the tautological bundle on the homogeneous space $G/P$ on which a maximal torus acts with characters $t_1,\dots, t_n$ and the $z_1, \dots ,z_n$ are formal complex variables. The Weyl group $W$ of $G$ acts on the set $\{ z_1, \dots , z_n \}$ by permutations, and the sets $X_i$ appearing in the formula are defined as $X_{i} = \{ \sigma(z_i): \sigma \in W \}$. The positive roots $\Phi^+ \setminus {\Phi_P}^+$ are expressed in the $z_i$ variables. \newline

In \cite{zielenkiewicz2014} we have obtained the above formulas in the case of classical Grassmannians, Lagrangian Grassmannians, orthogonal Grassmannians. The case of complete flag varieties is covered in \cite{berczi2012}. Recently we have obtained formulas for the partial flag varieties of types $A$, $B$, $C$ and $D$, using the methods presented here. The results are a part of the PhD dissertation of the author, submitted in December 2016. We also have computational results suggesting it holds for the homogeneous spaces $G_2/P_1$ and $G_2/P_2$, where $G_2$ is the smallest exceptional Lie group and $P_1, P_2$ are its two unique maximal parabolic subgroups. The aim of this paper is to describe an approach to the residue formulas which can hopefully be used to prove the above formula in general. \newline 

The meaning of the variables $z_i$ is not obvious from this description. Within this paper we show that they can be interpreted as the characters of an additional torus acting on $G$. The main idea behind this description is the Jeffrey--Kirwan nonabelian localization theorem. 

\subsection{Jeffrey--Kirwan nonabelian localization theorem}

Let $\X$ be a compact symplectic manifold equipped with a Hamiltonian action of a compact Lie group $G$, with  moment map $ \mu: \X \to \mathfrak{g}^* $. Assume $0$ is a regular value of $\mu$. One can form the symplectic reduction 
\[ \X \git G := \mu^{-1}(0) / G.\]
Then, there is a natural map $ \kappa: H^*_{G}(\X) \to H^*(\mu^{-1}(0)/G)$, called the Kirwan map, defined as the composition $(\pi^*)^{-1}\circ i^*$,
\[ \kappa: H^*_{G}(\X) \xrightarrow{i^*} H^*_{G}(\mu^{-1}(0)) \xrightarrow{(\pi^*)^{-1}}  H^*(\mu^{-1}(0)/G),\]
where $i^*$ is the map induced by the inclusion $i:\mu^{-1}(0) \to \X$ and $\pi^*$ is the natural isomorphism $H^*_{G}(\mu^{-1}(0)) \to H^*(\mu^{-1}(0)/G)$ induced by the quotient map $\pi: \mu^{-1}(0) \to \mu^{-1}(0)/G$. 
The assumption that $0$ is a regular value can of course be replaced by the assumption that $\xi \in \mathfrak{g}^*$ is a regular value and the reduction is taken at $\xi$, so that $\X \git_{\xi} G := \mu^{-1}(\xi) / G$. \newline

Assume $T$ is a maximal torus in $G$, and $\mu_T$ is the moment map for the $T$-action. 
The Jeffrey--Kirwan nonabelian localization theorem states that the Kirwan map is an epimorphism, and its kernel can be explicitly described in terms of intersection pairings. From the point of view of the residue formulas in equivariant cohomology, the most interesting result is the following (\cite{jeffrey1995}, Theorem 8.1):

\begin{theorem}[Jeffrey--Kirwan (\cite{jeffrey1995})] \label{thm:jk}

 Let $\omega$ be a symplectic form on $\X$ and $\omega_0$ the induced symplectic form on $\X \git G$. 
Let $\eta \in H^*_G(\X)$. Let $[\X \git G]$ be the fundamental class of $\X \git G$ in $H^*_G(\X)$. Let $\Phi^{+}$ and $W$ be, respectively, the set of positive roots and the Weyl group of $G$. Denote by $\varpi$ the product of the roots of $G$,\footnote{In the original formulation in \cite{jeffrey1995} $\varpi$ is the product of positive roots of $G$, so that the in the formulas $\varpi$ is replaced by $(-1)^{|\Phi^{+}|} \varpi^2$ .} $\varpi = \prod_{\gamma \in \Phi}\gamma$.
Then one can choose a subset $\mathcal{F}$ of the set of components of the fixed point set of the action of $T$ such that the following formula holds.
\[\kappa(\eta) e^{i \omega_0}[\X \git G] = \]
\[ = \frac{1}{(2 \pi )^{k-s} |W| \vol(T)} \Res\bigg{(}\varpi \sum_{F \in \mathcal{F}} e^{i \mu_T(F)} \int_F \frac{i^*_F(\eta e^{i \omega})}{e_F} \bigg{)}.\]
For $F \in \mathcal{F}$, the map $i_F$ is the inclusion of $F$ into $\X$ and $e_F$ is the equivariant Euler class of the normal bundle to $F$ in $\X$. The constant $\vol(T)$ is the volume of the torus $T$,  and $k,s$ denote the dimensions of $G$ and $T$ respectively.\footnote{The residue in the Jeffrey--Kirwan theorem is defined as a certain contour integral (see def. 8.5 in \cite{jeffrey1995}).} \newline

For a fixed component of the $T$ action, denote by $\beta_{F,j}$ the weights of the $T$-action on the normal bundle to $F$ in $\X$. Define $\beta^{\Lambda}_{F,j}$ to be equal $\beta_{F,j}$ if $\beta_{F,j} > 0$ and $-\beta_{F,j}$ otherwise. The set $\mathcal{F}$ of components of the fixed point set that appear in the summation is the set of those fixed components $F$ for which $\mu_T(F)$ lies in the cone $C_{F} := \{ \sum_{j} s_j \beta^{\Lambda}_{F,j}, s_j \geq 0\}$ spanned by the $\beta^{\Lambda}_{F,j}$.
\end{theorem}

One can hope to find a connection between the above theorem and the push-forward residue-type formulas via symplectic reductions. However, the Jeffrey--Kirwan localization theorem describes the non-equivariant cohomology of the symplectic reduction. Imposing an additional torus action on the space $\X$, and generalizing the Jeffrey--Kirwan theorem to the equivariant setting is the first step to obtain the residue-type push-forward formulas as special cases of the generalized Jeffrey--Kirwan formula. The preliminary results obtained in \cite{zielenkiewicz2014} in this setting suggest this approach will be successful. One can improve the Jeffrey--Kirwan theorem to work in $\T$-equivariant cohomology for an arbitrary torus $\T$ and one can reobtain the residue formula for classical Grassmannians from the Jeffrey--Kirwan theorem. \newline

In the assumptions of the Jeffrey--Kirwan theorem $X$ is a compact symplectic manifold. However, for our purposes it would be convenient to allow noncompact spaces. This can be done under the assumption that the moment map is proper. An explicit generalization of the Jeffrey--Kirwan theorem to noncompact spaces is presented in Kalkman's paper \cite{kalkman}, which uses the constructions in \cite{prato1994} to remove the compactess assumptions. A more detailed comment on that matter can be found in the end of Chapter \ref{gk_equiv} and an application is shown in Example \ref{grassmannian}.

\subsection{Other approaches to the nonabelian localization theorems}

An alternative statement and proof of the Jeffrey--Kirwan nonabelian localization theorem has been described by Guillemin and Kalkman in \cite{guillemin1996}. Their statement of the result is the following. Let $\X$ be a compact oriented symplectic manifold with an action of $T$, such that the $T$-action satisfies the assumptions of the Jeffrey--Kirwan localization theorem. Let 
\[\kappa: H^*_{T}(\X) \to H^*(\mu^{-1}(0)/T)\]
denote the Kirwan map.\newline

Consider the set  $\X_{crit} = \bigcup \X_j$ of critical points of the moment map $\mu$. Each $\X_j$ is a fixed point set of a one-dimensional subgroup $T_j$ of $T$. For a one-dimensional torus $T_j$ acting on $\X$, the equivariant cohomology can be computed from the Cartan complex
\[\tilde{\Omega} = \Omega^*(\X)^{T_j}\otimes \C[x_j],\]
where $\C[x_j]$ is a polynomial ring generated by an element $x_j$ of degree $2$, $x_j=c_1(\gamma_j)$ is the first Chern class of the universal complex line bundle $\gamma_j$ over $\C\P^{\infty} = \B T_j$. The differential in this complex is $\tilde{d} = d \otimes 1 + \iota(v)\otimes x_j$, where $d$ is the differential in $\Omega^*(\X)^{T_j}$ and $\iota(v)$ is the contraction with the vector field generating the $T_j$-action. Let $\kappa_j: H^*_{H_j}(\X_j) \to H^*(\X_j \git H_j)$ be the Kirwan map for the action of $H_j=T/T_j$ on $\X_j$. Denote by $i_j$ the inclusion $\X_j \to \X$, and let $e(\nu_j)$ be the Euler class of the normal bundle $\nu_j$ to $\X_j$ in $\X$. Define
\[\Res_j(\alpha):= \Res_{x_j =\infty}\frac{i^*_j \alpha}{e(\nu_j)},\]
where taking the residue at infinity means just taking the coefficient at $x_j^{-1}$ in the Laurent series (recall that we compute the equivariant cohomology from the Cartan complex, so we can express elements of $H^*_{H_j}(\X_j)$ as polynomials in $x_j$ with coefficients in $\Omega(\X_j)^{H_j}$).

\begin{theorem}[Guillemin--Kalkman (\cite{guillemin1996})]
\[\int_{\X \git T}\kappa(\alpha) = \sum_{j \in \mathcal{J}} \int_{\X_j \git H_j} \kappa_j (\Res_j(\alpha)).\]
The indexing set $\mathcal{J}$ is the set of those $j$ for which $X_j \cap \mu^{-1}(l) \neq \varnothing$ for a suitably chosen halfline $l \subseteq \mathfrak{t}^*$. For details see Section~\ref{gk_equiv}, equation~\eqref{jot}. 
\label{gkequiv}
\end{theorem} 
In the Guillemin--Kalkman theorem one only considers torus actions. If one considers an arbitrary compact group $G$, like in the Jeffrey--Kirwan theorem, then one can reduce the question of pushing-forward over $G$ to pushing-forward over its maximal torus $T$ using the Weyl Integration Formula. This is done in \cite{jeffrey1995} as a part of the proof of the Jeffrey--Kirwan theorem.

\begin{theorem}[Weyl Integration Formula \cite{weyl1925}]
Let $f$ be a continuous class function on a compact connected Lie group $G$ with maximal torus $T$ i.e. $f(t) = f(g t g^{-1})$. Then 
\[\int_{G} f(g) dg = \frac{(-1)^{n_+}}{|W|}\frac{vol(G)}{vol(T)}\int_{T} \varpi(t) f(t) dt,\]
where $\varpi(t) = \prod_{\gamma \in \Phi}\gamma(t), t \in T$ and $n_+$ is the number of positive roots of $G$.
\end{theorem} 

Using the above theorem, Jeffrey and Kirwan prove the following integration formula.

\begin{theorem}[Jeffrey--Kirwan (\cite{jeffrey1995})]
Let $\tilde{\alpha} \in H^*( \X \git T)$ be a lift of $\alpha \in H^*( \X \git G)$, in the sense that the pullback of $\tilde{\alpha}$ under the map induced by the canonical inclusion $i:\mu_{G}^{-1}(0)/T \to \X \git T$ equals the pullback of $\alpha$ under the fibration $p: \mu_{G}^{-1}(0)/T \to \X \git G$. Then
\[\int_{\X \git G} \alpha = \frac{1}{|W|} \int_{\X \git T} \tilde{\alpha} \cdot \varpi.\]
\label{MIT}
\end{theorem} 

Yet another viewpoint is presented by Martin in \cite{martin2000}.\footnote{The paper \cite{martin2000} was accepted for publication in Annals of Mathematics in December 1999 but has not appeared in print.} He proves that the cohomology of the symplectic reduction of a manifold $X$ by the action of $G$ can be described in terms of the cohomology of the symplectic reduction by the action of its maximal torus, as follows.

\begin{theorem}[Martin]
There is a natural ring isomorphism
\[H^*(\X \git G; \mathbb{Q}) \simeq \frac{H^*(\X \git T; \mathbb{Q})^W}{ann(\varpi)},\]
where $\varpi=\prod_{\alpha \in \Phi} \alpha$ and $ann(\varpi) = \{ c \in  H^*(\X \git T; \mathbb{Q})^W: c \cdot \varpi = 0\}$. \label{martin}
\end{theorem} 

Note that different authors use different conventions about measures on $G$ and $T$. The classical formulation of the Weyl Integration Formula uses explicitly the volumes of $G$ and $T$; Jeffrey and Kirwan assume the volume of $G$ is normalized to $1$, whereas Guillemin, Kalkman and Martin assume both the volumes of $G$ and $T$ to be normalized to $1$. 

\section{The equivariant Jeffrey--Kirwan localization theorem}\label{gk_equiv}

The Jeffrey--Kirwan nonabelian localization theorem can be improved to work in $\T$-equivariant cohomology. More precisely, let $\X$ be a compact oriented symplectic manifold equipped with a Hamiltonian action of $G$ with moment map $\mu_G$ and an action $\T$, such that the $G$-action satisfies the assumptions of the Jeffrey--Kirwan localization theorem and the actions of $G$ and $\T$ commute. We do not assume that the action of $\T$ is Hamiltonian. If the set $\mu_G^{-1}(0)$ is $\T$-invariant, then one can construct the equivariant Kirwan map
\[\kappa^{\T}: H^*_{G \times \T}(\X) \to H^*_{\T}(\mu_G^{-1}(0)/G),\]
defined as the composition 
\[H^*_{G \times \T}(\X) \xrightarrow{i^*} H^*_{G \times \T}(\mu_G^{-1}(0)) \xrightarrow{(\pi^*)^{-1}} H^*_{\T}(\mu_G^{-1}(0)/G), \]
where $i : \mu_G^{-1}(0) \to \X$ denotes the inclusion and $\pi:\mu_G^{-1}(0) \to \mu_G^{-1}(0)/G$ is the natural quotient map, inducing an isomorphism $ \pi^*: H^*_{G \times \T}(\mu_G^{-1}(0)) \to H^*_{\T}(\mu_G^{-1}(0)/G)$. In all of the above, $\mu_G$ denotes the moment map associated with the $G$-action. \newline

Some kind of equivariant Kirwan map has been introduced and investigated by Goldin in \cite{goldin2002}, in a very similar setting - in \cite{goldin2002} the equivariant Kirwan map is a map 
\[\kappa^{\T}: H^*_{K}(\X) \to H^*_{K/\T}(\mu^{-1}(0)/\T),\]
where $\T \lhd K$ is a subtorus of a compact Lie group $K$.  A proof of the surjectivity of the equivariant Kirwan map is given in \cite{goldin2002}. Its kernel is explicitly described, in terms analogous to Tolman-Weitsman's description of the kernel of the non-equivariant Kirwan map \cite{tolman2003}. The proof of the surjectivity presented in \cite{goldin2002} uses equivariant Morse theoretic methods to study the case of $\T= S^1$, followed by induction on the dimension of the torus. \newline

Theorem~\ref{thm:jk} gives a description of the non-equivariant push-forward. The effect of applying the push-forward map $\int_{\MG}:  H^*(\MG) \to H^*(pt)$ to the image of an element in $G$-equivariant cohomology under the Kirwan map $\kappa(\alpha)$ is described in terms of a certain residue operation. We relate this formula to the residue formulas obtained in \cite{zielenkiewicz2014}, by extending the results of Jeffrey and Kirwan to the equivariant setting. Our approach will be based on the paper by Guillemin and Kalkman \cite{guillemin1996}, in which the authors prove a formula similar to the one by Jeffrey--Kirwan. Their definition of a residue is slightly different, and they sum over a different subset of the fixed point set of the torus action. However, in the cases we will consider (projective spaces, Grassmannians) the two results coincide. We follow the notation of \cite{guillemin1996}, and we show how to adapt their proof of the residue-type push-forward formula in the equivariant setting. \newline

Let $\X$ be a compact symplectic manifold equipped with a Hamiltonian action of a compact group $G$ and an action of a torus $\T$, and assume the actions of $G$ and $\T$ commute and the set $\mu_G^{-1}(0)$ is $\T$-invariant. Assume $0$ is a regular value of the moment map $\mu_G $ and let $\MG := \mu_G^{-1}(0)/G $ denote the symplectic reduction of $\X$ for the $G$-action. Consider the $\T$-equivariant Kirwan map 
\[\kappa^{\T}: H^*_{G \times \T}(\X) \to H^*_{\T}(\MG),\]
and let $\int_{\MG}: H_{\T}^*(\MG) \to H_{\T}^*(pt)$ denote the equivariant push-forward. \newline

In view of the Jeffrey--Kirwan nonabelian localization theorem, it is possible to recover all results about $G$ from the action of the maximal torus $T$ in $G$. For this reason, let us assume that $G=T$, i.e. that $\X$ is equipped with two commuting actions of tori $T, \T$. The Borel model of the $\T$-equivariant cohomology is
 \[H^*_\T (\X) := H^*(\E\T \times^{\T} \X), \]
so to obtain the equivariant Guillemin--Kalkman theorem, one could try to apply the non-equivariant version of it to the space $\E\T \times^{\T} \X$. However, this space does not satisfy the compactness assumption. To avoid this problem, consider the finite-dimensional approximations $\E_m$ such that 
\[ H^i_G(\X) = H^i(\E_m \times^G \X) \textrm{ for } m \textrm{ sufficiently large with respect to }i. \]
Typically, for the torus actions one takes $\E_m = (\C^m \setminus 0)^n$ with the action of $\T$ given by multiplication on every component. However, these spaces are not compact, so it is better to consider spheres $\E_m = (S^{2m+1})^n$. \newline

The following $\T$-equivariant version of the Guillemin--Kalkman formula holds. Let $\X_{crit} = \bigcup \X_j$ be the set of critical points of the moment map $\mu_T$. Each $\X_j$ is a fixed point set of a one-dimensional subgroup $T_j$ of $T$. Let $\kappa^{\T}_j$ be the equivariant Kirwan map for the action of $H_j=T/T_j$ on $\X_j$,
\[\kappa^{\T}_j: H^*_{H_j \times \T}(\X_j) \to H^*_{\T}(\X_j \git H_j),\]
and let $\kappa^\T$ be the $\T$-equivariant Kirwan map for the action of $T$ on $\X$,
\[\kappa^{\T}: H^*_{T \times \T}(\X) \to H^*_{\T}(\X \git T).\] 
Denote by $i_j$ the inclusion $X_j \to X$ and let $e^{\T}(\nu_j)$ be the equivariant Euler class of the normal bundle $\nu_j$ to $\X_j$ in $\X$. Let $x_j$ be a chosen generator of $H^*_{T_j}(pt)$. Define
\[\Res_j(\alpha):= \Res_{x_j=\infty}\frac{i^*_j \alpha}{e^{\T}(\nu_j)}.\] 

\begin{theorem}[Equivariant Guillemin--Kalkman Theorem for $S^1$-actions] With the assumptions and notation introduced above, the following formula holds.
\[\int_{\X \git T}\kappa^{\T}(\alpha) = \sum_{j \in \mathcal{J}} \int_{\X_j \git H_j} \kappa^{\T}_j (\Res_j(\alpha)).\]
The indexing set $\mathcal{J}$ is the set of those $j$ for which $X_j \cap \mu^{-1}(l) \neq \varnothing$ for a suitably chosen halfline $l \subseteq \mathfrak{t}^*$. For details see Section~\ref{gk_equiv}, equation~\eqref{jot}. 
\label{gkequiv}
\end{theorem}

The proof is a straightforward application of the Guillemin--Kalkman result to the space $\tilde{\X} = \E_m \times^{\T} \X$ for a sufficiently large $m$ (which depends on $\alpha$), where $\E_m = (S^{2m+1})^n$. Let us choose $\alpha \in H_{\T}^k(\X)$ and choose $m$ such that $H_{\T}^k(\X) \simeq H^k(\E_m \times^{\T} \X)$, identifying $\alpha$ with an element in $H^k(\E_m \times^{\T} \X)$. Then the push-forward 
\[H^*_{\T}(\X \git T) \to H^*_{\T}(pt)=H^*(\B\T)\]
can be described as follows. The equivariant cohomology of a point $H^*_{\T}(pt)$ is a polynomial ring $\C[t_1, \dots, t_n]$. Let us denote $t=(t_1,\dots, t_n)$ and let $I=(i_1, \dots, i_n)$ be a multi-index such that $|I| = \frac{1}{2}(\deg \alpha - \dim \X \git T )$. Then the push-forward sends
\[\kappa^{\T}(\alpha) \mapsto \sum a_I t^I,\]
and the coefficients $a_I$ are given by the integrals
\[\int_{\X_I } \kappa^{\T}(\alpha),\]
where $\X_I = (S^{2i_1+1} \times \dots \times S^{2 i_n +1}) \times^{\T} \X \git T$.  \newline

The proof of Guillemin and Kalkman is based on induction. The $S^1$ case is described in detail in the Appendix (adapting the Guillemin--Kalkman argument to the equivariant case). For one-dimensional torus actions the argument is more general: it does not require symplecticity assumptions, one only needs a compact orientable manifold with boundary $(\X, \partial \X)$ together with an $S^1$-action which is locally free on the boundary. Instead of the symplectic reduction one can then take $\partial \X / S^1$. However, to proceed with the induction one needs to choose subsequent one-dimensional tori in $T$ in such a way that in every step the assumptions are satisfied, i.e. one needs to choose a splitting $T=S \times H$, where $S$ is a one-dimensional torus, acting locally freely on $\partial \X / H$. If we knew that $\partial \X / H$ is the boundary of a compact manifold $\partial \X / H = \partial M$, then using the Guillemin--Kalkman theorem for $S=S^1$, we could express the integral 
\[\int_{\partial \X / S \times H} \kappa(\alpha) = \int_{\partial M / S} \kappa(\alpha) = \sum_k \int_{M_k} \Res_{x_k=\infty} \frac{i^*_k \alpha}{e(\nu_k)},\]
where the $M_k$ are the connected components of the fixed point set of $S^1$. In general one cannot claim that $\partial \X / H $ is the boundary of some compact manifold. This is where we need the symplectic structure and the moment map of the action. \newline

The assumptions that $\X$ is symplectic and the action of $T$ is Hamiltonian enable one to use the moment map for the action for choosing the one-dimensional subtori. Using the moment map, one can give a precise description of how to proceed with induction. The argument is based on the following theorem, due to Atiyah \cite{atiyah1982} and independently to Guillemin and Sternberg \cite{guillemin1982}. If $\X$ is a compact manifold equipped with a Hamiltonian torus action with moment map $\mu$ and $0$ is a regular value of $\mu$, then the image of the moment map is a convex polytope $\Delta \subseteq \lie{t}^*$. We assume $0$ lies in its interior, otherwise the symplectic reduction is trivial, i.e. $\X \git T=\varnothing$ and Theorem \ref{gkequiv} is tautologically true. The set $\Delta^0$ of the regular values of $\mu$ is a disjoint union of convex polytopes 
\[\Delta^0 = \Delta^0_1 \cup \dots \cup \Delta^0_k,\]
and by assumption $0$ lies inside one of the $\Delta^0_j$. For any chosen element $\theta$ in the weight lattice of $\lie{t}^*$, one can consider a ray through the origin in the direction of $\theta$,
\begin{equation}\label{jot}
l = \{ s \theta: s \in [0, \infty) \}.
\end{equation}
Let us choose $\theta$ in such a way that the ray $l$ does not intersect any of the walls of $\Delta_i^0$ of codimension greater than one and hence intersects the codimension one walls transversely. Then the Lie subalgebra $\lie{h} \subseteq \lie{t}$ defined as
\[\lie{h} = \{ v \in \lie{t}: \langle \theta, v \rangle = 0 \}\]
is the Lie algebra of a codimension one subtorus $H \subseteq T$. The assumptions we made on the ray $l$ and hence on the element $\theta$ imply that the moment map $\mu$ is transverse to $l$ and the action of $H$ on $\mu^{-1}(l)$ is locally free. Moreover, the action of $T/H \simeq S^1$ on $\mu^{-1}(l)/H$ is locally free on the boundary, which makes it possible to proceed with induction. Note that $\mu^{-1}(l)/H$ might not be a symplectic manifold, but rather a symplectic orbifold. However, the proof can easily be adapted to the case of orbifolds, because most of the analysis is done locally, and the only global component of the proof is Stokes Theorem, which holds for orbifolds \cite{satake1957}. The details of the proof for the $S^1$ case are described in the Appendix. \newline

Following the procedure which is briefly described above and in full detail in \cite{guillemin1996}, one arrives at the following result. Let us choose the ray $l$ as described above, and call $l$ a main branch. Next, for every intersection point $p_j$ of $l$ with the codimension one walls of $\Delta_i^0$, let us choose a ray $l_j$ starting at $p_j$ and not intersecting any codimension three walls of $\Delta_i^0$. The rays $l_j$ are called secondary branches. Next, continue the procedure by considering the intersection points of secondary branches $l_j$ with codimension 2 walls of $\Delta_i^0$, and at each such point choose ternary branches (rays not crossing the codimension 4 walls), etc. Finally one arrives at a vertex of the moment polytope (the vertices of the moment polytope are the images of fixed components). One obtains what Guillemin and Kalkman in \cite{guillemin1996} called a dendrite $\mathcal{D}$ - a set of branches, consisting of sequences of rays $(l, l^{(1)}, \dots , l^{(n)})$, where $l^{(j)}$ is a branch constructed in step $j$, and a set of points $(0, p_1, \dots, p_n)$ on those branches, such that the branch $l^{(i)}$ starts at the point $p_i$ and intersects the codimension $i+2$ wall at $p_{i}$. Each branch $B$ in the dendrite $\mathcal{D}$ corresponds to a fixed component of the action of $T$ and determines a sequence of tori
\[\{0\} \subseteq T_F^{(1)} \subseteq T_F^{(2)} \subseteq \dots \subseteq T_F^{(n)} = T,\] 
which gives the choice of the desired one-dimensional subtori $T_F^{(i+1)}/T_F^{(i)}$ needed for the induction. Given a sequence of tori $T_F^{(i)}$ one can choose a basis $\mathbf{x}_F= \{x_{F,1},\dots,x_{F,n}\}$ of $\lie{t}^*$ such that for each $i$ the dual elements $x_{F,1}^*, \dots, x_{F,i}^*$ form a basis of the integer lattice in the Lie algebra $\lie{t}_F^{(i)}$ of the torus $T_F^{(i)}$.

\begin{theorem}[Equivariant Guillemin--Kalkman Theorem] With the notation introduced above, if one additionally assumes that the fixed points of the action are isolated\footnote{If the fixed points are not isolated one needs to replace the summation over fixed points with summation over the fixed components, followed by the push-forward to a point over the component.} and denote by $i_p$ the inclusion of the fixed point $p$ into $\X$ and by $e(\nu_p)$ the Euler class of the normal bundle at $p$, applying theorem \ref{MIT} gives:

\[\int_{\X \git G}\kappa^{\T}(\alpha) = \frac{1}{|W|}\sum_{B \in \mathcal{D}} \Res_{\mathbf{x}_p=\infty}\bigg{(}\varpi \frac{i^*_p \alpha}{e(\nu_p)}\bigg{)}.\]

\end{theorem}

\begin{remark}
The dendrite procedure, introduced by Guillemin and Kalkman in \cite{guillemin1996}, is generalized by Jeffrey and Kogan in \cite{jeffrey2005} by replacing rays with higher dimensional cones, yielding a different proof of the Jeffrey--Kirwan nonabelian localization theorem.
\end{remark}

\begin{example}[Projective plane]
Consider $\C^2$ with a $T = S^1$-action given by $x(z_0, z_1) = (x z_0, x z_1)$ and an $\T = S^1 \times S^1$-action given by $(t_0, t_1)(z_0, z_1) = (t_0 z_0, t_1 z_1)$. The moment map to the $S^1$-action is 
\[ \mu (z_0, z_1) = (|z_0|^2 + |z_1|^2), \]
so $\mu^{-1}(1) = S^3 \subseteq \C^2$ and $\mu^{-1}(1)/S^1 = \C\P^1$.\newline

Let $\X = \{ z \in \C^2: \mu(z) \leq 1 \}$ be a ball of radius $1$ in $\C^2$, so that $\X \git S^1 = \partial \X / S^1 = \C\P^1$. Let $\alpha \in H^k_{S^1 \times \T}(\X)$ and consider 
\[\int_{\C\P^1} \kappa^\T(\alpha).\]
The $S^1 \times \T$-equivariant cohomology of $\X$ is a module over $H^*_{S^1 \times \T}(pt) = \C[x,t_0,t_1]$.
The only fixed point for the action of $S^1$ on $\X$ is the origin, so the expression
\[\sum_{\{  k : \mu_{|\X_k} \geq 0 \}} \int_{\X_k} \Res_{x=\infty} \frac{i^*_k \alpha}{e(\nu_k)}\]
reduces to a sum over a single point $0=(0,0)$ and the integral is just the evaluation at this point. We get:
\[\int_{\C\P^1} \kappa^\T(\alpha) = \Res_{x=\infty} \frac{i^*_0 \alpha}{e(\nu_0)}.\]\newline
The restriction of the form $\alpha \in H^k_{S^1 \times T}(\X)$ to a point $0$ is a polynomial in $x$, let us denote this polynomial by $f$. The $S^1 \times \T$-equivariant Euler class equals $(t_0-x)(t_1 -x)$. One has:
\[\int_{\C\P^1} \kappa^\T(\alpha) = \Res_{x=\infty} \frac{f(x)}{(t_0 - x) (t_1 - x)}=\]
\[ = \frac{f(t_0)}{ (t_1 - t_0)} + \frac{f(t_1)}{(t_0 - t_1)} = \frac{f(t_0)-f(t_1)}{ (t_1 - t_0)} ,\]
which is a polynomial in $t_0, t_1$. \newline

Alternatively, one could consider the action of $S^1$ on $\C\P^2$, via
\[z[z_0:z_1:z_2]  = [z_0, z z_1, z z_2],\]
and use dendrites to describe the push-forward. Then the fixed points of the $S^1$ action are $[1:0:0]$ and $\{ [0:z_1:z_2] \} = \C\P^1$ and only the fixed point $[1:0:0]$ appears in the summation.
\end{example}

\begin{example}[Grassmannian]
Consider the action of unitary matrices $U(k)$ on the space $Hom(\C^k,\C^n)$ given by multiplication on the right. The moment map for this action \linebreak $\mu: Hom(\C^k, \C^n) \to \mathfrak{u}(k)$ is given by
\[ \mu(A) = A^* A - Id, \]
hence $\mu(0)^{-1} = \{ A^*A = Id\}$ and the column vectors of a matrix $A \in \mu(0)^{-1}$ form a unitary $k$-tuple in $\C^n$, so
\[Hom(\C^k, \C^n) \git U(k) = Grass_k(\C^n).\]
Let $T$ denote the maximal torus in $U(k)$, acting on $Grass_k(\C^n)$ via restriction of the action of $U(k)$, with characters $\mathbf{z} = \{ z_1,\dots,z_k\}$. The roots of $U(k)$ are $ \Phi = \{ z_i - z_j\}_{i \neq j}$. Applying Theorem \ref{MIT} we can reduce the integral over $Grass_k(\C^n)$ to the integral over the reduction of $Hom(\C^k, \C^n)$ with respect to the action of $T$:

\[ \int_{Grass_k(\C^n)} \kappa_\T(\alpha) = \frac{1}{|W|} \int_{Hom(\C^k, \C^n) \git T} \varpi  \kappa^T_\T(\alpha) , \]
where $\varpi = \prod_{\gamma \in \Phi}\gamma = \prod_{i \neq j} (z_i - z_j)$ is the product of the roots of $U(k)$, $W$ denotes the Weyl group of $U(k)$ and $\kappa^T_\T$ is the $\T$-equivariant Kirwan map for the action of $T$ on $Hom(\C^n, \C^k)$. \newline

The maximal torus $T$ in $U(k)$ acts diagonally on $Hom(\C^k, \C^n) = \C^n \oplus \dots \oplus \C^n$, with the action on each component $\C^k$ given by multiplication. The moment map $\mu_T : Hom(\C^k, \C^n) \to \mathfrak{t}^*$ for this action is given by projection of the moment map for the action of $U(k)$, hence
\[\mu_T(A) = (||v_1||,\dots, ||v_k||),\]
where $v_1, \dots, v_k$ are the column vectors of $A$. The symplectic reduction for the $T$-action is therefore
\[Hom(\C^k, \C^n) \git T = \mu_T^{-1}(0)/T = (\C\P^{n-1})^k,\]
since $\mu_T^{-1}(0)$ consists of $k$-tuples of vectors of length $1$ in $\C^n$. The image of the moment map $\mu_T$ is the product of half-lines $(\RR_{\geq 0})^k$. It is not a convex polytope, because the space $Hom(\C^k, \C^n)$ is non-compact. However, we can still use the dendrite algorithm as in \cite{guillemin1996}, only this time we have to make sure that at each step we choose rays $l$ in such a way, that they intersect a codimension one wall (and not diverge to infinity). In our case this can be easily achieved - if the chosen ray $l$ does not intersect a codimension one face of $(\RR_{\geq 0})^k$, then the ray $-l$ does. Choosing the rays in this way will always lead to a branch ending at the only fixed point of the action, the origin. Therefore, one gets the following expression for the push-forward

\[ \int_{Grass_k(\C^n)} \kappa_\T(\alpha) = \frac{1}{|W|}\Res_{\mathbf{z}=\infty} \frac{\varpi V(z_1,\dots, z_k)}{e(0)}, \]
where $\varpi  = \prod_{i \neq j} (z_i - z_j)$ is the product of the roots of $U(k)$, and $e(0)$ is the $\T \times T$-equivariant Euler class at zero, $e(0) = \prod_{i,j}(z_i - t_j)$, where $z_1,\dots,z_k$ are the characters of the $T$-action and $t_1,\dots,t_k$ are the characters of the $\T$-action. Finally,

\[ \int_{Grass_k(\C^n)} \kappa_\T(\alpha) = \frac{1}{|W|} \Res_{\mathbf{z}=\infty} \frac{ \prod_{i \neq j} (z_i - z_j) V(z_1,\dots,
z_k)}{\prod_{i,j}(z_i - t_j)}, \]
which is the same formula as the one obtained in \cite{zielenkiewicz2014}, using purely combinatorial methods. \newline

The description of the Grassmannian as the symplectic reduction for the $U(k)$-action has been used by Martin to give a description of the non-equivariant cohomology of the classical Grassmannian and to derive an integration formula using Theorem \ref{martin}. For details, see \cite{martin2000}, chapter 7.
\label{grassmannian}
\end{example} 

\section{Appendix}\label{appendix}

\subsection{The outline of the proof}

Our aim is to apply the Jeffrey--Kirwan nonabelian localization theorem (in the formulation of Guillemin and Kalkman) to reobtain the residue-type formulas for push-forwards in equivariant cohomology of classical Grassmannians. The key idea is to realize the Grassmannian as a symplectic reduction. If $\X$ is a symplectic manifold with a Hamiltonian action of a compact group $G$, then the Jeffrey--Kirwan theorem provides a description of the non-equivariant push forward $H^*(\X \git G) \to H^*(pt)$ in terms of a certain residue. \newline

We assume that $\X$ is a symplectic manifold with action of a product $G \times \T$, where $G$ is a compact group with maximal torus $T$ and acts in a Hamiltonian way, and $\T$ is a compact torus. We assume that $\mu^{-1}_G(0)$ and $\mu_{T}^{-1}(0)$ are $\T$-invariant. We rephrase the nonabelian localization theorem in $\T$-equivariant cohomology (adapting the proof of Guillemin and Kalkman), hence obtaining an expression for the equivariant push-forward $H^*_{\T}(\X \git G) \to H^*_{\T}(pt)$ in terms of a certain residue. In the case of the classical Grassmannian, the obtained formula coincides with the one obtained in \cite{zielenkiewicz2014}. \newline

The proof we present here uses two main tools:

\begin{enumerate}
	\item Theorem \ref{MIT}, which reduces the problem to the case $G=T$.
	\item The non-equivariant approach of Guillemin and Kalkman applied to the approximation spaces of the Borel model of the $\T$-equivariant cohomology of $\X$.
\end{enumerate}
Guillemin and Kalkman study the Kirwan map by looking at it on the level of differential forms. They give a precise description of the isomorphism 
\[ \pi^*: H^*(\partial \X / T) \to H^*_{T}(\partial X),\]
in terms of iterated residues, by first describing $\pi^*$ for $T=S^1$ and using induction on the dimension of the torus. Finally, they apply the Stokes theorem and the Atiyah--Bott--Beline--Vergne localization formula to obtain the formula for the push-forward. \newline

In the equivariant case we need to describe the isomorphism
\[ \pi^*_{\T}: H^*_{\T}(\partial \X / T) \to H^*_{T \times \T}(\partial \X),\]
and use it to deduce the formula for the $\T$-equivariant push-forward \linebreak $H^*_{\T}(\partial \X ) \to H^*_{\T}(pt)$. For this we consider the approximation spaces of the Borel model of the $\T$-equivariant cohomology of $\X$, $X_m = \E_m \times^{\T} \X$ with $E_m = (S^{2m + 1})^n$ and use the following fact relating the equivariant push-forward with the non-equivariant one.

\begin{proposition}
Let $\X$ is a compact manifold with a $\T$-action. Let $p: \X \to pt$, and let $\alpha \in H^k_{\T}(\X)$. Then the push-forward $p_*: H^*_\T(\X) \to H^*_\T(pt)$ is given by the formula
\[p_* \alpha = \sum_I \beta_I t^I,\]
where $I = (i_1, \dots, i_n)$ is the multi-index satisfying $|I| = \frac{1}{2}(\deg \alpha- \dim \X)$. The coefficients $\beta_I$ are given by 
\[\beta_I = \int_{\X_I} j^* \alpha,\]
with $\X_I = ( S^{2i_1 + 1} \times S^{2 i_2 + 1} \times \dots \times S^{2 i_l + 1} ) \times^\T \X \xrightarrow{j} \E\T \times^\T \X$ .
\end{proposition}

\subsection{The $S^1$ action - a detailed proof via approximation spaces.}
 
Let $\X$ be a compact oriented $S^1$-manifold with boundary, and assume that the $S^1$ action is locally free on the boundary of $\X$.  The inclusion $ i: \partial{\X} \to \X$ induces the map on $S^1$-equivariant cohomology \[i^*: H^*_{S^1}(\X) \to  H^*_{S^1}(\partial{\X}).\]
 Recall that if the action on $\partial{\X}$ is locally free, there is a canonical isomorphism $\pi^*: H^*(\partial{\X} / S^1 ) \to H^*_{S^1}(\partial{\X})$. The Kirwan map is defined as the composition 
\[ \kappa = (\pi^*)^{-1} \circ i^* : H^*_{S^1}(\X) \to H^*(\partial{\X} / S^1 ). \]
Equivariant cohomology can be computed from the Cartan complex, which in the case of an $S^1$-action on $\partial{\X}$ is 
\[ \tilde{\Omega} = \Omega^*(\partial{\X})^{S^1} \otimes \C[x] \]
with differential $\tilde{d} = d \otimes 1 + \iota(v) \otimes x$, where and $\iota(v)$ denotes the contraction with the vector field generating the $S^1$-action. \newline

Guillemin and Kalkman in \cite{guillemin1996} consider the integral $\int_{\partial{\X} / S^1} \kappa(\alpha)$ for an equivariantly closed form $\alpha \in \tilde{\Omega}^k$. Since our aim is to prove an analogous result for push-forwards in equivariant cohomology, we replace the Kirwan map $\kappa$ by its equivariant analogue $\kappa_{\T}: H^*_{S^1 \times \T} \X \to H_{\T}^* \partial{\X} / S^1 $ and consider the push-forward in $\T$-equivariant cohomology
\[\int_{\partial{\X} / S^1} \kappa_{\T}(\alpha),\]
for an $(S^1 \times \T)$-equivariantly closed form $\alpha$. We will always assume that $\T$ is a torus and that the actions of the two tori $S^1$ and $\T$ commute. The proof presented here follows the proof in \cite{guillemin1996}, adapting it to the equivariant setting. \newline

If $\X$ is a compact $S^1 \times \T$-manifold with a locally free $ S^1$ action on $\partial{\X}$ and the actions of $S^1$ and $\T$ commute, then the spaces $\E_m \times^{\T} \X$ are compact $S^1$-manifolds with a locally free action of $S^1$ on the boundary. Moreover, $\partial{\E_m \times^{\T} \X} = \E_m \times^{\T} \partial \X$ are the approximation spaces in the Borel model for equivariant cohomology of $\partial \X$. Let us denote $\X_m =  \E_m \times^{\T} \X$ and let $\tilde{\Omega}_{\partial \X_m}$ be the Cartan complex computing the $S^1$-equivariant cohomology of $\partial \X_m$ (this implies, that for $m$ large enough it computes the $S^1 \times \T$-equivariant cohomology of $\partial \X$ in degrees $i$ small with respect to $m$). Let $\alpha \in \tilde{\Omega}_{\partial \X_m}^*$ be an equivariantly closed form, with respect to the $S^1$-action, and assume $\deg \alpha = \dim \partial \X_m -1 $. Then, using the isomorphism
\[ \pi^*: H^*(\partial \X_m / S^1) \simeq H^*_{S^1}(\partial \X_m),\]
we can write 
\[\alpha = \tilde{d} \nu + \pi^* \gamma\]
for some $\nu \in \tilde{\Omega}^{k-1}_{\partial \X_m}$ and $\gamma \in \Omega^k(\partial \X_m / S^1)$.
Following the calculation in section 2 of \cite{guillemin1996} we can explicitly write down $\nu$ and $\pi^* \gamma$.
Consider the following element of $\Omega^*(\partial \X_m)^{S^1}\otimes \C[x,x^{-1}]$ (which is the just the localization of $\tilde{\Omega}_{\partial \X_m}$ with respect to $x$):
\[ \nu_0 = \frac{\theta \alpha }{x + d \theta} = \frac{\theta \alpha}{x} \sum_{n \geq 0} \big{(}\frac{-d \theta}{x}\big{)}^n \]
where $\theta$ is an $S^1$-invariant 1-form on $\partial \X_m$ satisfying $\iota(v)\theta = 1$.
The element $\nu_0$ is a Laurent series in $x$ with coefficients in $\Omega^*(\partial \X_m)^{S^1}$. Writing the form $\alpha = \sum_{i \geq 0} \alpha_i x^i$ as a polynomial in $x$ with coefficients in $\Omega^*(\partial \X_m)^{S^1}$, we can rewrite $\nu_0$ as
\[ \nu_0 = \sum_{n,i} \theta \alpha_i (-d \theta)^n x^{i-n-1}.\]
Since the coefficients of the above series are in $\Omega^*(\partial \X_m)^{S^1}$, which is trivial in degrees higher than $\dim \partial \X_m$, the coefficient of $x^{i-n-1}$ is zero if
\[1 + \deg \alpha_i + 2n > \dim \partial \X_m. \]
Since $\deg \alpha = \dim \partial \X_m -1$, the only non-zero coefficients appear when $n \leq i$, so the only powers of $x$ that can occur in the Laurent series of $\nu_0$ are the positive ones and $x^{-1}$, so we can write
\[\nu_0 = \nu + \beta x^{-1}, \]
where $\beta = \Res_{x=\infty} \nu_0 \in \Omega^*(\partial \X_m)^{S^1}$. \newline

The differential in the Cartan complex $\tilde{\Omega}_{\partial \X_m} = \Omega^*(\partial \X_m)^{S^1} \otimes \C[x]$ is given by $\tilde{d} = d \otimes 1 + \iota(v)\otimes x$, so 
\[\alpha = \tilde{d} \nu + \iota(v)\beta,\]
and the form $\iota(v)\beta$ is $S^1$-invariant and horizontal, so it comes from some form $\gamma \in \Omega^k(\partial \X_m / S^1)$,
\[\iota(v)\beta = \pi^* \gamma.\]
This shows that the map $(\pi^*)^{-1}$ is given by the formula  
\[(\pi^*)^{-1}(\alpha) = \Res_{x=\infty} \iota(v) \frac{\theta \alpha}{x + d \theta} = \Res_{x=\infty} \pi_{*}  \frac{\theta \alpha}{x + d \theta}.\]
 \newline

Let $\{ (\X_m)_i \}_{i=1, \dots, N}$ be the connected components of the fixed point set of the action of $S^1$, and let $\{ U_i \}_{i=1, \dots, N}$ be pairwise disjoint tubular neighbourhoods of sets $ (\X_m)_i$ such that $U_i \cap \partial \X_m = \varnothing$. Let $\alpha \in H^*_{S^1}(\X)$ (we will abuse notation and denote by $\alpha$ also its restriction to $X_m$ and $\partial \X_m$), and let $\nu,\theta$ be as above. Assume that $\deg \alpha = \dim \partial \X_m -1$. Let us extend the forms $ \nu , \theta$ to $\X_m \setminus  \X_m^{S^1}$. Applying the Stokes theorem to $ 0 = \int_{\X_m} \alpha =  \int_{\X_m} \tilde{d} \nu $ we get
\[ \sum_{k=1}^N \int_{U_k} \frac{\theta \alpha}{x + d \theta} = \int_{\partial \X_m} \frac{\theta \alpha}{x + d \theta} = \int_{\partial \X_m / S^1} \pi_{*}\big( \frac{\theta \alpha}{x + d \theta} \big).\]
It was shown in \cite{berline1983}\footnote{The computation is a step of the proof of a theorem announced in \cite{berline1982} and proven in \cite{berline1983}. The most detailed computation can be found in \cite{guillemin1999supersymmetry}.} that by shrinking the radii of $U_i$ to zero the left-hand side converges to
\[ \sum_{k=1}^N \int_{(\X_m)_k} \frac{i^*_k \alpha}{e(\nu_k)},\]
where $i_k : (\X_m)_k \to \X_m$ is the inclusion map and $e(\nu_k)$ denotes the Euler class of the normal bundle to $(\X_m)_k$ in $\X_m$. Taking residues at $x=\infty $ of both sides of the above expression we get
\[ \int_{\partial \X_m / S^1} \Res_{x=\infty} \pi_{*}\big( \frac{\theta \alpha}{x + d \theta} \big) = \sum_{k=1}^N \int_{(\X_m)_k} res_{x=\infty} \frac{i^*_k \alpha}{e(\nu_k)}, \]
and the right hand-side equals $\int_{\partial \X_m / S^1} \kappa^\T (\alpha)$. \newline

Finally, we have shown that for $\alpha \in H^{\dim \partial \X_m - 1}(\partial \X_m)$ we have
\[\int_{\partial \X_m / S^1} \kappa^{\T} (\alpha)  = \sum_{k=1}^N \int_{(\X_m)_k} \Res_{x=\infty} \frac{i^*_k \alpha}{e(\nu_k)},\]
so for a class $\alpha \in H^{i}_T(\partial \X)$ we can choose $m$ such that $i=\dim \X_m -1$ and it follows that for  $\alpha \in H^i_{\T}(\partial \X)$ the equivariant push-forward satisfies 
\[\int_{\partial \X/ S^1} \kappa^{\T} (\alpha)  = \sum_{k=1}^N \int_{\X_k} \Res_{x=\infty} \frac{i^*_k \alpha}{e(\nu_k)}.\]

Note that our notation for the residues is different than the one of Guillemin and Kalkman, who denote the residues by $\Res_{x=0}$, not by $\Res_{x=\infty}$. However, this is a matter of notation only. The residue in \cite{guillemin1996}, like here, is defined to be the coefficient at $x^{-1}$ in the series expansion, which we choose to call the residue at infinity to remain consistent with the classical notation from calculus. 

\begin{corollary}[Symplectic case]
Consider a symplectic manifold $\X$ with a Hamiltonian action of $S^1 \times \T$ with moment map for the $S^1$-action $\mu: \X \to \lie{s}^*$. Then the manifold with boundary $\X_+ = \{ x \in \X : \mu(x) \geq 0 \}$ is a compact $S^1$ manifold with a locally free $S^1$-action on the boundary $\partial \X_+ = \mu^{-1}(0)$. Let $\X \git S^1$ denote the quotient $\mu^{-1}(0)/ S^1$ and let $\X_k$ denote the connected components of the fixed point set of the action on $\X \git S^1$. By the considerations above applied to $\X_+$ we get:

\[\int_{\X \git S^1} \kappa^{\T}(\alpha) = \sum_{\{ k: \mu_{|\X_k} \geq 0 \}} \int_{\X_k} \Res_{x=\infty} \frac{i^*_k \alpha}{e(\nu_k)}. \] 
\end{corollary}

\begin{corollary}[Jeffrey--Kirwan formula]
In \cite{jeffrey1995} the authors consider the case when $\X$ is a symplectic manifold with an action of a compact group $K$ with the moment map $\mu: M \to \lie{k}^*$. The differential form which is being integrated is
\[\kappa(\alpha) = \eta_0 e^{i \omega_0},\]
where $\eta_0 = \kappa(\eta)$  is the image of the Kirwan map for some form $\eta \in H^*_{K}(\X)$ and $\omega_0$ is the symplectic form on the symplectic reduction $\X \git K$. In the special case when $K = \T = S^1$, one obtains the following formula for the push-forward:
\[\int_{\X \git S^1} \eta_0 e^{i \omega_0} = \sum_{F \subseteq F_+} \int_{F} \Res_{x=\infty} \frac{i^*_F (\eta e^{i (\omega + \mu)})}{e(\nu_F)} =  \]
\[ = \Res_{x=\infty} \sum_{F \subseteq F_+} \int_{F} e^{i \mu(F) }\frac{i^*_F (\eta e^{i \omega})}{e(\nu_F)},\]
where $F \subseteq F_+$ are the components of the subset $F_+$ of the fixed-point set on which $\mu$ is positive. \newline
 
After a slight change of notation, this is almost the formula from Theorem 8.1 in \cite{jeffrey1995} and its corollary below, corrected by eliminating the $\frac{1}{2}$-factor and the $\psi^2$-factor under the residue. Again, we denote the series coefficient at $x^{-1}$ as a residue at infinity, not at zero as in \cite{jeffrey1995}.
\end{corollary}

\bibliographystyle{alpha}
\bibliography{bibliography}

\end{document}